\documentclass[oneside,a4paper,titlepage]{article}
\usepackage[T1]{fontenc}
\usepackage[latin1]{inputenc}
\usepackage{amsmath}
\usepackage{amssymb}
\usepackage{amsthm}

\newtheorem{ex}{Example}

\newtheorem{Th}{Theorem}
\newtheorem{lem}{Lemma}

\newcommand{\RR}{\mathbb{R}}

\newcommand{\NN}{\mathbb{N}}

\newcommand{\cD}{{\cal D}}

\newcommand{{{\cadlag}}}{c\`adl\`ag}

\begin{document}

\sffamily
\pagestyle{plain}

\title{A brief note on the soundness \\of Bermudan option pricing via cubature}
\author{F S Herzberg, Merton College, University of Oxford, Oxford OX1 4JD}

\maketitle

\begin{abstract} 
The subject of this study is an iterative Bermudan option pricing algorithm based on (high-dimensional) cubature. We show that the sequence of Bermudan prices (as functions of the underlying assets' logarithmic start prices) resulting from the iteration is bounded and increases monotonely to the approximate perpetual Bermudan option price; the convergence is linear in the supremum norm with the discount factor being the convergence factor. Furthermore, we prove a characterisation of this approximated perpetual Bermudan price as the smallest fixed point of the iteration procedure.
\end{abstract}

\noindent

When Nicolas Victoir studied ``asymmetric cubature formulae with few points'' \cite{V} for symmetric measures such as the Gaussian measure, the idea of (non-perpetual) Bermudan option pricing via cubature in the log-price space was born. In the following, we will discuss the soundness and convergence rate of this approach when used to price perpetual Bermudan options.

Consider a convex combination $(\alpha_1,\dots,\alpha_m)\in[0,1]^d$ (that is, $\sum_{k=1}^m \alpha_k = 1$) and $x_1,\dots,x_m\in\RR^d$. Then there is a canonical {\em weighted arithmetic average operator} $A$ associated with $\vec \alpha, \vec x$ given by $$\forall f\in\RR^\RR\quad Af=\sum_{k=1}^m\alpha_k f(\cdot-x_k).$$ Now suppose $c\in(0,1)$, $g,h:\RR^d\rightarrow \RR$, $Ag\geq g$, $Ah=h$ and $0\vee g\leq h$. Define an operator $\cD$ on the cone of nonnegative measurable functions by $$\cD:f\mapsto \left(c\cdot Af\right)\vee g.$$ Since $A$ is positive and linear, thus monotone (in the sense that for all $f_0\leq f_1$, $A f_0\leq A f_1$), it follows that $\cD$ must be monotone as well. Furthermore, whenever $Af\geq f$, we have that $A\cD f\geq \cD f$, as the linearity and positivity of $A$ combined with our assumption on $g$ imply $$A\cD f \geq cA^2f\vee g\geq \left(c\cdot Af\right)\vee g =\cD f.$$ Finally, due to our assumptions on $h$ and $g$, we have for all nonnegative $f\leq h$, $$\cD f\leq cAh\leq h.$$ Summarising this, we are entitled to state

\begin{lem} \label{discreteLemma}Adopting the previous paragraph's notation and setting $$Q:=\left\{f\leq h \ : \ Af\geq h\right\},$$ we have that $$\cD:Q\rightarrow Q,$$ $\cD$ is monotone (i e order-preserving), and $A\cD -\cD$ is nonnegative.
\end{lem}

This is sufficient to prove 

\begin{Th} \label{monoconv}For all $n\in\NN_0$, \begin{equation}\label{discretemono}\cD^{n+1}(g\vee 0)\geq \cD^n(g\vee 0)=:q_n.\end{equation} Furthermore, $$q:=\lim_{n\rightarrow \infty}\cD^n(g\vee 0)=\sup_{n\in\NN_0}\cD^n(g\vee 0)\in Q$$ and $q$ is the smallest nonnegative fixed point of $\cD$.
\end{Th}
\begin{proof}
\begin{enumerate}
\item The proof of equation (\ref{discretemono}) is a straightforward induction on $n$ where we have to use the monotonicity of $\cD$ in the induction step. 
\item Since $\cD$ maps $Q$ itself, the whole sequence $\left({\cD}^{n}(g\vee 0)\right)_{n\in\NN_0}$ is bounded by $h$. This entails $q\leq h$ as well. Using the linearity of $\sup$ and our previous observation that $A\cD-\cD\geq 0$ (Lemma \ref{discreteLemma}), we can show \begin{eqnarray*}\forall y\in\RR^d\quad Aq(y)&= &\sup_{n\in\NN_0} A\left({\cD}^{n}(g\vee 0)\right)(y)\\ &\geq& \sup_{n\in\NN_0}{\cD}^{n}(g\vee 0)(y) =q(y),\end{eqnarray*} which means $Aq\geq g$. As we have already seen, $q\leq h$, so $q\in Q$.
\item Again, due to the linearity of $\sup$, $\cD$ and $\sup_{n\in\NN_0}$ commute for bounded monotonely increasing sequences of functions. Thereby $$\cD q=\sup_{n\in\NN_0}\cD{\cD}^{n}(g\vee 0)=\sup_{n\in\NN} {\cD}^{n}(g\vee 0)=q.$$
\item Any nonnegative fixed point $p$ of $\cD$ must be greater or equal $g\vee 0$. Therefore by the monotonicity of $\sup$ and $\cD$, $$\sup_{n\in\NN_0}{\cD}^{n}p\geq \sup_{n\in\NN_0}{\cD}^{n}(g\vee 0)=q.$$
\end{enumerate}
\end{proof}

\begin{lem} \label{gLemma} Using the previous Theorem's notation, we have for all $x\in\RR^d$ and $n\in\NN_0$, if $q_{n+1}(x)=g(x)$, then $q_n(x)=g(x)$.
\end{lem}
\begin{proof} By the monotonicity of the sequence $(q_n)_{n\in\NN_0}$ (Theorem \ref{monoconv}), we have $$g(x)\leq q_0(x)\leq q_n(x)\leq q_{n+1}(x).$$
\end{proof}

\begin{Th} For all $n\in\NN$, $$\left\|q_{n+1}-q_n\right\|_{C^0\left(\RR^d,\RR\right)}\leq c\cdot \left\|q_{n}-q_{n-1}\right\|_{C^0\left(\RR^d,\RR\right)}.$$
\end{Th}

\begin{proof} The preceding Lemma \ref{gLemma} yields \begin{eqnarray*} \left\|q_{n+1}-q_n\right\|_{C^0\left(\RR^d,\RR\right)}&=& \left\|q_{n+1}-q_n\right\|_{\left\{q_{n+1}>g\right\}}\\ &=&\left\|c\cdot Aq_{n}-\left(\left(c\cdot A q_{n-1}\right)\vee g\right)\right\|_{C^0\left(\left\{c\cdot Aq_{n}>g\right\},\RR\right)}\end{eqnarray*} via the definition of $q_{i+1}$ as $\left(cAq_i\right)\vee g$ for $i=n$ and $i=n+1$. But the last equality implies \begin{eqnarray*}\left\|q_{n+1}-q_n\right\|_{C^0\left(\RR^d,\RR\right)}&\leq&\left\|c\cdot Aq_{n}-c\cdot A q_{n-1}\right\|_{C^0\left(\left\{c\cdot Aq_{n}>g\right\},\RR\right)}\\ &\leq&\left\|c\cdot Aq_{n}-c\cdot A q_{n-1}\right\|_{C^0\left(\RR^d,\RR\right)}.\end{eqnarray*} Since $A$ is linear as well as an $L^\infty$-contraction (and therefore a $C^0$-contraction, too), we finally obtain $$ \left\|q_{n+1}-q_n\right\|_{C^0\left(\RR^d,\RR\right)}\leq c\left\|A\left(q_{n}-q_{n-1}\right)\right\|_{C^0\left(\RR^d,\RR\right)}\leq c\left\|q_{n}-q_{n-1}\right\|_{C^0\left(\RR^d,\RR\right)}.$$
\end{proof}

\begin{ex}[Bermudan put option with equidistant exercise times in $t\cdot \NN_0$ on the weighted arithmetic average of a basket in a discrete Markov model with a discount factor $c=e^{-rt}$ for $r>0$] Let $\beta_1,\dots,\beta_d\in[0,1]$ be a convex combination and assume that $A$ is such that \begin{equation}\label{discretecondition}\forall i\in \{1,\dots,d\}\quad \sum_{k=1}^m\alpha_ke^{-(x_k)_i}=1,\end{equation} then the functions $$g:x\mapsto K-\sum_{i=1}^d\beta_i \exp\left(x_i\right)$$ and $h:=K$ (where $K\geq 0$) satisfy the equations $Ah=h$ and $Ag=g$, respectively. Moreover, by definition $g\leq h$. Then we know that the (perpetual) Bermudan option pricing algorithm that iteratively applies $\cD$ to the payoff function $g\vee 0$ on the $\log$-price space, will increase monotonely and will have a limit which is the smallest nonnegative fixed point of $\cD$. Moreover, the convergence is linear and the contraction rate can be bounded by $c$.

The condition (\ref{discretecondition}) can be achieved by a change of the time scale (which ultimately leads to different cubature points for the distribution of the asset price)
\end{ex}

{\bf Acknowledgements.} This work originates from research conducted by the author for his doctoral thesis at the University of Oxford. The author gratefully acknowledges funding from the German Academic Exchange Service ({\em Doktorandenstipendium des Deutschen Akademischen Austauschdienstes}) and helpful discussions with Professor Terry Lyons.

\end{document}